\documentclass[a4paper,12pt]{article}
\usepackage{amssymb}
\usepackage{amsthm}
\theoremstyle{plain}
\newtheorem{Lemma}{Lemma}
\newtheorem{Proposition}{Proposition}
\newtheorem{Theorem}{Theorem}

\theoremstyle{definition}

\newtheorem{Example}{Example}

\newcommand{\car}{\mathrm{char}}
\newcommand{\PP}{\mathbb{P}}

\newcommand{\K}{\mathbb{K}}

\begin{document}
\title{Birational geometry of defective varieties}
\author{Edoardo Ballico and Claudio Fontanari}
\date{}
\maketitle 

\begin{small}
\begin{center}
\textbf{Abstract}
\end{center}
Here we investigate the birational geometry of projective varieties of 
arbitrary dimension having defective higher secant varieties.
We apply the classical tool of tangential projections and we determine 
natural conditions for uniruledness, rational connectivity, and 
rationality. 
\vspace{0.5cm}


\end{small}

\vspace{0.5cm}

\section{Introduction}

We work over an algebraically closed field $\K$ with $\car(\K)=0$. 
Let $X \subset \PP^r$ be an integral nondegenerate $n$-dimensional 
variety. Recall that for every integer $k \ge 0$ the 
$k$-\emph{secant variety} $S^k(X)$ is defined as the Zariski closure 
of the set of the points in $\PP^r$ lying in the span of $k+1$ 
independent points of $X$. An easy parameter count shows that the 
expected dimension of $S^k(X)$ is exactly $\min \{r, n(k+1)+k \}$. 
However, there are natural examples of projective varieties having 
secant varieties of strictly lower dimension: for instance, the 
first secant variety of the $2$-Veronese embedding of $\PP^2$ in 
$\PP^5$ has dimension $4$ instead of $5$. More generally, one 
defines the $k$-\emph{defect} $\delta_k(X) = \min \{r, n(k+1)+k \} 
- \dim S^k(X)$ and says that $X$ is $k$-\emph{defective} 
if $\delta_k(X) \ge 1$. It seems reasonable to regard defective 
varieties as exceptional and try to classify them. 
The first result in this direction, stated by Del Pezzo in 1887 
and proved by Severi in 1901 (see \cite{DelPezzo:87}, \cite{Severi:01}, 
and also \cite{Dale:85} for a modern proof), characterizes the $2$-Veronese 
of $\PP^2$ as the unique $1$-defective surface which is not a cone. Along 
the same lines, subsequent contributions by Palatini (\cite{Palatini:06}
and \cite{Palatini:09}), Scorza (\cite{Scorza:08} and \cite{Scorza:09}),
and Terracini (\cite{Terracini:11} and \cite{Terracini:21}) set up 
the classification of $k$-defective surfaces and of $1$-defective 
varieties in dimension up to four. This classical work has been 
recently reconsidered and generalized by Chiantini and Ciliberto 
in \cite{CC1} and \cite{CC2}. 
It turns out that one of the main tools for understanding defective 
varieties is provided by the technique of tangential projections. 
Namely, assume that $X$ is not $(k-1)$-defective and let $p_1, \ldots, 
p_k$ be general points on $X \subset \PP^r$. The general 
$k$-\emph{tangential projection} $\tau_{X,k}$ is the projection  
of $X$ from the linear span of its tangent spaces at $p_1, \ldots, 
p_k$. By the classical Lemma of Terracini (see \cite{Terracini:11} 
and \cite{Dale:84} for a modern version), $X_k := \tau_{X,k}(X)$ 
is lower dimensional than $X$ if and only if $X$ is $k$-defective. 
Therefore, the classification of defective varieties reduces to the 
classification of varieties which drop dimension in the general 
tangential projection. However, we believe that the only reasonable 
goal in arbitrary dimension is to determine some geometrical properties 
of defective varieties, and the present paper 
is indeed a first attempt in this direction. In order to state our 
main results, we recall from \cite{CC1} that the \emph{contact locus} 
of a general hyperplane section $H$ tangent at $k+1$ general points 
$p_1, \ldots, p_{k+1}$ of $X$ is the union $\Sigma$ of the irreducible 
components of the singular locus of $H$ containing $p_1, \ldots, p_{k+1}$. 
One defines $\nu_k(X) := \dim \Sigma$ and says that $X$ is $k$-\emph{weakly 
defective} whenever $\nu_k(X) \ge 1$. The reason for this terminology 
is simply that a $k$-defective variety is always weakly defective 
(indeed, we are going to show in Proposition~\ref{defects} that 
$\nu_k(X) \ge \delta_k(X)$), but the converse is not true (look at 
cones). We point out the following: 

\begin{Lemma}\label{birational} 
Fix integers $k \ge 1$, $n \ge 2$, $r \ge (k+1)(n+1)$, and let 
$X \subset \PP^r$ be an integral nondegenerate $n$-dimensional variety. 
If $X$ is not $k$-weakly defective, then $\tau_{X,k}$ is birational.
\end{Lemma}

Moreover, for every $k \ge 1$ and $n \ge 2$ we exhibit a projective 
$n$-dimensional variety being $k$-weakly defective but not $k$-defective 
such that  $\tau_{X,k}$ is not birational (see Example~\ref{counter1}). 
Next, as an application of Lemma~\ref{birational}, we investigate the 
birational geometry of a $k$-defective variety of arbitrary dimension: 

\begin{Theorem}\label{uniruled}
Fix integers $k \ge 1$, $n \ge 2$, $r \ge (k+1)(n+1)-1$, and let 
$X \subset \PP^r$ be an integral nondegenerate $n$-dimensional variety
which is $k$-defective but not $(k-1)$-defective.
Assume that $\nu_k(X) = \delta_k(X)$ and, for $k \ge 2$, that 
$X$ is not $(k-1)$-weakly defective. Then $X$ is uniruled. 
Assume moreover that the general contact locus of $X$ is irreducible. 
Then $X$ is rationally connected and for $\nu_k(X) = \delta_k(X)=1$ 
it is rational.
\end{Theorem}

We stress that the hypotheses for uniruledness cannot be removed: 
in Examples~\ref{counter2} and \ref{counter3} we collect a series 
of non-uniruled defective varieties of any dimension which do 
not satisfy exactly one of the listed assumptions. An analogous remark 
applies to the second part of the statement, where irreducibility 
turns out to be essential. For instance, if $X$ is a cone 
over a curve of positive genus, then two general points on $X$ cannot be 
joined by a rational curve. Indeed, we suspect that a defective variety 
with $\nu_k = \delta_k$ and reducible general contact locus should be a 
cone (this is certainly the case in dimension up to four by \cite{CC2}, 
proof of Proposition~4.2, and \cite{Scorza:09}, \S~14).

This research is part of the T.A.S.C.A. project of I.N.d.A.M., 
supported by P.A.T. (Trento) and M.I.U.R. (Italy).

\section{The proofs}

Let $k$ be the minimal integer such that $X$ is $k$-defective. 
Then $\tau_{X,k-1}(X)$ is generically finite and by definition 
we have 
\begin{equation}\label{one}
\nu_k(X)=\nu_1(X_{k-1}).
\end{equation}
Moreover, by applying twice the equality $\delta_h(Y) = \dim Y 
- \dim \tau_{Y,h}(Y)$, first with $Y = X$ and $h=k$, then with 
$Y = X_{k-1}$ and $h=1$, we deduce
\begin{equation}\label{two}
\delta_k(X)=\delta_1(X_{k-1}).
\end{equation}

\begin{Proposition}~\label{defects}
Fix integers $k \ge 1$, $n \ge 2$, $r \ge (k+1)(n+1)-1$, and let 
$X \subset \PP^r$ be an integral nondegenerate $n$-dimensional variety.
Assume that $k$ is the minimal integer such that $X$ is $k$-defective. 
Then we have $\nu_k(X) \ge \delta_k(X)$.
\end{Proposition}

\proof Assume first $k=1$. If $\delta_1(X)=1$, we are just claiming 
that $\nu_1(X)>0$, which is well-known (see for instance \cite{CC1}, 
Theorem~1.1). If instead $\delta_1(X)\ge 2$, we take a general hyperplane 
section $H$. By \cite{CC2}, Lemma~3.6, we have  $\nu_1(H)=\nu_1(X)-1$ 
and $\delta_1(H)=\delta_1(X)-1$, so we conclude by induction. 
Assume now $k > 1$. By (\ref{one}) and (\ref{two}), we are reduced to 
the previous case, so the proof is over.

\qed

\emph{Proof of Lemma~\ref{birational}.} Recall that $\tau_{X,k}$ is
defined as the projection of $X$ from $<~T_{p_1}(X), \ldots, T_{p_k}(X)~>$, 
where $p_i$ is a general 
point on $X$ and $T_{p_i}(X)$ denotes the tangent space to $X$ at $p_i$.   
Pick a general point $p_{k+1} \in X$ and let $q = \tau_{X,k}(p_{k+1})$. 
By Proposition~\ref{defects}, $X$ is not $k$-defective, so $\tau_{X,k}$ 
is generically finite. Therefore $\tau_{X,k}^{-1}(q)= \{ p_{k+1}, \ldots, 
p_{k+d} \}$ and we want to show that $d=1$. 
Since $\tau_{X,k}(T_{P_{k+h}}(X)) = T_q(X_k)$ 
for every $h$ with $1 \le h \le d$, it follows that 
$< T_{p_1}(X), \ldots, T_{p_k}(X), T_{P_{k+h}}(X) >$ does not 
depend on $h$. On the other hand, by \cite{CC1}, Theorem~1.4, 
the general hyperplane which is tangent to $X$ at $p_1, \ldots, 
p_{k+1}$ is not tangent to $X$ at any other point, so we have 
$d = h =1$ and the proof is over. 

\qed

In order to construct nontrivial examples of projective varieties 
whose general $k$-tangential projection is not birational, we are 
going to apply the following criterion:

\begin{Lemma}~\label{criterion} 
Fix integers $k \ge 1$, $n \ge 2$, $r \ge (k+1)(n+1)-1$, and let 
$X \subset \PP^r$ be an integral nondegenerate $n$-dimensional variety.
Assume that $X$ is $k$-weakly defective, but not $k$-defective. 
If $X$ is not uniruled, then $\tau_{X,k}(X)$ is generically finite 
but not birational.
\end{Lemma}

\proof Notice that $X$ is a fortiori not $(k-1)$-defective, so we can 
apply Proposition~3.6 in \cite{CC1} with $h=k$ to obtain that $X_k$ is 
$0$-defective. Hence $X_k$ is a developable scroll (see for instance 
\cite{CC1}, Remark~3.1~(ii)), in particular it is uniruled and it 
follows that $X_k$ is not birational to $X$. 

\qed

\begin{Example}~\label{counter1} Let $k \ge 1$, $n \ge 2$, and 
$r \ge (k+1)(n+1)-1$. Take a $(n-1)$-dimensional variety $C$ in 
$\PP^r$ which is not uniruled, a linear space $V \subset \PP^r$ 
with $\dim V = k$, and a smooth hypersurface $H \subset \PP^r$ of degree 
$d \ge k+2$ such that $V \not\subset H$. Let $W$ be the cone over 
$C$ with vertex $V$ and define $X := H \cap W$. By \cite{CC1}, 
Example~4.3, $X$ is $k$-weakly defective but not $k$-defective. 
We claim that $X$ is not uniruled. Let $\pi: X \to C$ be the 
projection and take a general point $p \in X$. If $R \subseteq X$ 
is a rational curve passing through $p$, then $R$ is contained in 
a fiber of $\pi$, otherwise its projection would be a rational 
curve through a general point of $C$, in contradiction with the 
non-uniruledness of $C$. On the other hand, the general fiber of 
$\pi$ is set-theoretically the intersection of $H$ with a 
$(k+1)$-dimensional linear space, hence it is a smooth hypersurface 
of high degree, which is not covered by rational curves. Hence the 
claim is established and from Lemma~\ref{criterion} we deduce that  
$\tau_{X,k}(X)$ is generically finite but not birational.
\end{Example}

\emph{Proof of Theorem~\ref{uniruled}.} Assume first $k = 1$. 
If $\nu_1(X) = \delta_1(X) = 1$, from the so-called Terracini's 
Theorem (\cite{CC1}, Theorem~1.1) it follows that the general contact locus 
$\Sigma$ of $X$ imposes only three conditions on hyperplanes containing 
it, in particular it is a plane curve. Moreover, if $\Sigma$ had degree 
$d > 2$, then the general secant line to $X$ would be a multisecant, 
which is a contradiction (the above argument is borrowed from 
\cite{CC2}, proof of Proposition~4.2). 
Hence $X$ is uniruled and if $\Sigma$ is irreducible then $X$ is 
rationally connected. More precisely, if $d = 1$ then two 
general points on $X$ would be joined by a straight line, a 
contradiction since $X$ is nondegenerate. Therefore if we fix a 
general $p \in X$ and for general $q \in X$ we assume that 
the corresponding contact locus $\Sigma_{pq}$ is irreducible, 
then $\Sigma_{pq}$ is a smooth conic and a natural birational map 
from $X$ to its tangent space $T_p(X) \cong \PP^n$ is defined by 
sending $q$ to the intersection point between the tangent lines 
to $\Sigma_{pq}$ at $p$ and $q$ (see \cite{Scorza:09}, \S~15). 
If instead $\delta_1(X) \ge 2$, let $H$ be a general hyperplane 
section. By \cite{CC2}, Lemma~3.6, we have  $\nu_1(H)=\nu_1(X)-1$ 
and $\delta_1(H)=\delta_1(X)-1$, so uniruledness and rational 
connectivity follow by induction. Assume now $k > 1$. 
By (\ref{one}) and (\ref{two}), the previous cases apply to 
$X_{k-1}$. On the other hand, from Lemma~\ref{birational} we 
get a birational map between $X$ and $X_{k-1}$, so the proof is over.

\qed 

\begin{Example}\label{counter2} Here we show that the assumption $X$ 
not $(k-1)$-weakly defective is essential for every $n \ge 2$. 
Fix $r \ge 3n+2$ and let $C \subset \PP^r$ be a $(n-1)$-dimensional 
variety which is not uniruled. Take a line $L \subset \PP^r$ 
and a smooth hypersurface $H \subset \PP^r$ of degree $d \ge 3$ 
such that $L \not\subset H$. Let $W$ be the cone over 
$C$ with vertex $V$ and define $X := H \cap W$. By \cite{CC1}, 
Example~4.3, we have $\delta_2(X) = \nu_2(X) = 1$ and $X$ is also 
$1$-weakly defective. Moreover, arguing as in Example~\ref{counter1}, 
it is easy to check that $X$ is not uniruled. 
\end{Example}

Finally we focus on the assumption $\nu_k(X) = \delta_k(X)$. 
By Proposition~\ref{defects}, it is always satisfied in the 
case of surfaces; however, in higher dimension it is no more 
automatic. 

\begin{Example}\label{counter3} As in \cite{CC2}, Example~2.2, 
we consider a $n$-dimensional variety $X$ contained in a 
$(n+1)$-dimensional cone $W$ over a curve $C$ in $\PP^r$ 
with $r \ge 2n+1$. We have $\delta_1(X)=n-2$ and 
$\nu_1(X)=n-1$; moreover, if we assume $g(C)\ge 1$ and we 
take $X := H \cap W$, where $H$ is a general hypersurface 
in $\PP^r$ of high degree, then the same argument as in 
Example~\ref{counter1} shows that $X$ is not uniruled.
\end{Example}

\noindent
Edoardo Ballico \newline
Universit\`a degli Studi di Trento \newline
Dipartimento di Matematica \newline
Via Sommarive 14 \newline
38050 Povo, Trento, Italy \newline
e-mail: ballico@science.unitn.it

\vspace{0.5cm}

\noindent
Claudio Fontanari \newline
Universit\`a degli Studi di Trento \newline
Dipartimento di Matematica \newline
Via Sommarive 14 \newline
38050 Povo, Trento, Italy \newline
e-mail: fontanar@science.unitn.it

\end{document}